\magnification=\magstep1
\baselineskip=18pt
\font\bsc=cmcsc10 scaled\magstep2

\def\proclaim #1. #2\par{\medbreak
\noindent{\bf#1. \enspace}{\sl#2}\par\medbreak}

\def\section #1. #2\par{\bigbreak
\noindent{\bf#1. \enspace}{\bf#2}\par\medbreak} 

\def\C{{\cal C}}

\def\vo{v^{(1)}}

\def\del{\partial}
\def\v2{v^{(2)}}
\def\bv2{\bar{v}}
\def\bbv2{\hat{v}}
\def\t2{t^{(2)}}
\def\bt2{\bar{t}}
\def\bbt2{\hat{t}}
\def\ga2{\gamma}
\def\pA{A^{\prime}}
\def\g{\gamma}
\def\R{{\bf R}}
\def\df{{\dot f}}
\def\ddf{{\ddot f}}

\def\tt{{\tilde t}}
\def\tv{{\tilde v}}

\def\pt{t^{\prime}}
\def\pv{v^{\prime}}
\def\py{y^{\prime}}
\def\r{\rightarrow}
\def\e{\epsilon}
\def\N{{\bf N}}

\def\qed{\vrule height 5pt width 5 pt depth 0pt}

\centerline{\bsc acceleration of bouncing balls in external fields}
\medskip
\centerline{T.~Kr\"uger,\footnote{$^1$}{Forschungszentrum BiBoS, 
Universit\"at Bielefeld}
L.D.~Pustyl'nikov$^{1,}$\footnote{$^2$}{Electric Power Research Institute 
VNIIE, Moscow} and S.E.~Troubetzkoy$^{1,}$\footnote{$^3$}{Institute for 
Mathematical Sciences, SUNY Stony Brook}}
\medskip
\baselineskip=13pt
\midinsert
\narrower
We introduce two models, the Fermi-Ulam model in an external field and a one
dimensional system of bouncing 
balls in an external field above a periodically oscillating
plate.  For both models we investigate the possibility of unbounded motion.
In a special case the two models are equivalent.

AMS Classification: 70
\endinsert
\baselineskip=17pt
\section 1. Introduction
\par
In this article we
investigate the possibility of unbounded growth of energy in periodically
forced Hamiltonian systems. The first models of this type were proposed by 
Poincar\'e [P] in relation to the problem of the growth of entropy. In [Pu8] 
these questions were studied in Newtonian as well as relativistic
mechanics.
One of the most popular simple models for which this question has been 
investigate is the Fermi-Ulam model [F][U]: a point particle moves vertically
between two periodically oscillating plates and reflects elastically from
them.  For rationally related frequencies of the oscillating plates it was shown
[Pu4][Pu5][Pu9] that the velocity of the ball is bounded for each trajectory.  
Another popular model is that of an elastic ball bouncing vertically on a
periodically
oscillating plate under the influence of an external gravitational field.  In
contrast to the Fermi-Ulam model it was shown in [Pu1]-[Pu3] that for
a sufficiently fast oscillating plate for an open set of initial condition with
infinite measure the velocity of the ball is unbounded.

In this paper we propose two new models for the investigation of Fermi 
acceleration which are natural generalizations of the above described 
models.  The first is the Fermi-Ulam model in an external gravitational
field which will be described precisely
in section 2.  The second is a one dimensional 
system of bouncing balls in an external 
gravitational field above a periodically oscillating plate which will be 
described precisely in section 3. For each of these models we investigate the 
question, does there exist a trajectory with unbounded growth of 
the absolute value of the velocity and the energy. 
Furthermore we show that in a certain special case the Fermi-Ulam model in 
an external field is equivalent to the one dimensional system of two balls.
In contrast to the Fermi-Ulam model, for the bouncing ball model we study
only special cases.

\section 2. The Fermi-Ulam model in an external field with and without 
singularities
\par
At first we describe the (classical) Fermi-Ulam model [U]. 
We consider two infinitely heavy parallel plates, the lower plate 
oscillates according to the law $z=f_1(t)$ and the upper plate according to the 
law $z=f_2(t).$  Here $t$ is time and we assume the both functions are 
smooth, periodic with period $T > 0$ and $f_2(t) > f_1(t).$ Note that the period
of both functions is assumed the same! A ball moves freely
in the direction $z$ between the plates and collides elastically with them.
Here the laws of motion and collision are Newtonian.  The Newtonian 
mechanics play an important role, in [Pu6][Pu7] a relativistic 
analogue of the Fermi-Ulam model was introduced and it was shown that 
the energy of the ball has unbounded growth. Below we will describe the 
behavior of the classical Fermi-Ulam model which is in stark contrast to the
relativistic case.

Ulam [U] posed the following problem in connection with mechanism of Fermi 
acceleration [F]: do there exist trajectories for which 
the energy or the velocity of the ball grow unboundedly.  This 
problem for the Fermi-Ulam model was first solved affirmatively in [Pu4][Pu5]
for analytic functions $f_1(t)$ and $f_2(t).$ Ulam himself only considered 
the case when $f_1(t)$ and $f_2(t)$ are sinusoidal and moreover the special 
case when one of them is constant.  The main result of [Pu4][Pu5] is:
for any initial condition the energy and the velocity 
of the ball are bounded by a certain constant which depends on the initial 
conditions. The proofs in [Pu4][Pu5] actually hold for
smooth enough differentiable functions $f_1(t)$ and $f_2(t).$

In this article we study the Fermi-Ulam model in an external gravitational 
field, the ball falls toward the lower plate with constant 
acceleration $g \ge 0$ (fig.~1). 
The main result for this model is that if $f_1(t)$ and $f_2(t)$ are smooth 
enough and do not touch (i.e. $f_1(t) < f_2(t)$ for all $t$) then the 
energy and the velocity of the ball are always 
bounded, i.e.~are smaller than a constant which does not depend on $t.$

We say that the orbit of the ball which has speed $v_0$ at time $t_0$ has 
{\it bounded velocity} if there exists $c := c(v_0) > 0$ such that 
$\sup_{\tau} |v(\tau)| \le c.$ If the ball does not have bounded velocity we
say that it has {\it unbounded velocity}.
For simplicity we state the next theorem for analytic functions even though 
the proof holds for finitely differentiable functions with enough derivatives.

\proclaim Theorem 2,1. Suppose $f_1(t), f_2(t)$ are analytic, periodic with 
the same period and $f_1(t) < f_2(t)$ for all $t.$ Then the ball has 
bounded velocity.\par

Since the ball has zero radius we consider an important special case when the 
two plates are allowed to touch which we call a singularity, i.e.~there
exists a time $t^*$ when $f_1(t^*) = f_2(t^*)$ (fig.~2).
In this case the ball can make an 
infinite number of collisions with both plates in a finite time interval.
We say that the ball which has speed $v$ at time $t$ has {\it unbounded 
velocity at time} $t^* > t$ if $\lim_{\tau \nearrow t^*} |v(\tau)| 
= \infty.$

\proclaim Theorem 2,2.
Suppose $f_1,f_2$ are continuous function and there is an $\epsilon_1 > 0$
such that for $k=1$ or $k=2$ the $f_i$ are $C^k$ functions on $(t^*
-\epsilon_1,t^*).$ Furthermore assume that $f_1(t) < f_2(t)$ for $t
\in (t^* - \epsilon_1,t^*),$ and $f_i,f_i^{(1)},\dots,f_i^{(k)}$
extend by continuity from the left to the points $t^*.$ If
$f^{(j)}_1(t^*) = f^{(j)}_2(t^*)$ for $j=0,1 \dots k-1$ and $f_1^{(k)}(t^*)
\ne f_2^{(k)}(t^*)$ then
there exist positive constants $K,\epsilon$ such that if the 
initial condition $(v_0,t_0)$ satisfies $|v_0| > K$ and 
$t_0 \in (t^*-\epsilon,t^*)$ 
then the particle will have unbounded velocity at time $t^*$. 
\par
\noindent
Theorems 2,1 and 2,2 will be proven in section 4.
\section 3. A one dimensional system of bouncing balls in an external
gravitational field with oscillating plates 
\par
In this article we will only consider the case of two balls.  Call the 
balls  $P_1$ and $P_2$ and the masses $m_1$ 
and $m_2.$  The balls are elastic point particles and move on a vertical 
line in a gravitational field with constant acceleration $g > 0.$  They 
collide elastically with each other and with an infinitely heavy horizontal 
plate which oscillated periodically in time according to the law $z = f(t).$
Here $f(t)$ is a smooth function having period $T > 0.$ 
In the case of triple collisions, standard mechanical laws 
do not apply in most cases.  
In the special case when instead of two ball we have 
one ball then this model reduces to the model studied in [Pu1]-[Pu3].
We remark that several balls under the influence of an external field 
bouncing vertically above a fixed plate was studied in [W1][W2] and [Ch].
Our first result is
\proclaim Proposition 3,1. If the masses of both balls are positive but finite 
then either both balls have unbounded velocity or both balls have bounded
velocity.\par

To produce acceleration we consider a special class of functions which was
introduced in [Pu1]-[Pu3].
Let $\C := \{ f(t) : f \hbox{ periodic with period } T \hbox{ and } \exists
t_0$ with $\df(t_0) = KTg/2$ for some $K \in \N^+\}.$ Note that
if $f$ is of class $C^1$ then there exists a positive constant $c_0(f)$
such that $f(ct) \in \C$ for all $c \le c_0$ and $f(ct) \not \in \C$ for all
$c > c_0,$ that is if we scale $f$ to oscillate quickly then it is in the class
$\C$ and if we scale $f$ to oscillate slowly then it is not in $\C.$ The
constant $c_0(f)$ depends only of the variation of $f.$

We next consider the special case when the balls are of
equal masses.  We investigate which different kinds of possible acceleration can
happen. Namely suppose that $P_2$ is the lower ball and $P_1$ is the upper ball.
Let $\v2_i $ be the velocity of ball $P_2$ at the moment just after the
$i$th collision with the plate and
$\vo_i $ be the velocity of the ball $P_1$ at the moment just after the
$i$th collision with ball $P_2.$
\proclaim Theorem 3,2. If the masses of both balls are equal and if
$f \in \C$ then there exists unbounded motions of both
balls such that
\item{(a)} $\lim_{i \r \infty} |\vo_i| = \lim_{i \r \infty} |\v2_i| = \infty,$
\item{(b)} $\limsup_{i \r \infty} |\vo_i|= 
\limsup_{i \r \infty} |\v2_i| = \infty,$
\item{} $\max (\liminf_{i \r \infty}|\vo_i|,
\liminf_{i \r \infty}|\v2_i|) < \infty,$\par

Next we consider the special cases when the mass of ball $P_1$ 
is zero while $P_2$ has non-zero mass. In analogy to 
celestial mechanics this case is called the restricted system.
There are two possible positions of the balls, in the first case the ball 
$P_2$ with non-zero mass is situated between the plate 
and the ball $P_1$ (fig.~3) and in the second case the ball $P_1$ is situated 
between the plate and the ball $P_2$ (fig.~4). In both case the two balls 
interact between each other according to the law of elastic collisions.
The ball $P_1$ does not influence the ball $P_2$ since it has zero mass.

\proclaim Theorem 3,3. If $f$ is a periodic function of class $C^1$ then 
in case (1) it is impossible that $P_2$ has unbounded 
velocity while $P_1$ has bounded velocity. 
If $f \in \C$ then there are initial conditions for which both balls have 
unbounded velocity.\par
\noindent
If $f \not \in \C$ does the motion (unbounded, unbounded) occur?
Arnold has shown for a class ${\cal D} \subset \C$ of slowly oscillating
functions the ball $P_2$ always has bounded velocity [A].
\proclaim Theorem 3,4. 
In case (2), if $f$ is a periodic function of class $C^1$ then
there exists bounded motions of $P_2$ and if $f \in \C$ then there exists 
unbounded motions of $P_2,$ for each of  
which an open set with infinite measure of initial conditions for ball $P_1$
lead to $P_1$ having infinite velocity in finite time.\par
Note that the motion of ball $P_2$ can, in this case, be regularized through
a triple collision since $P_1$ has zero mass and does not affect $P_2.$
The motion of $P_1$ can not be regularized through this collision.
In case (2) if $P_2$ has periodic motion then the model is equivalent to a 
special case of the Fermi-Ulam model in an external gravitational field.
\section 4. Proofs of Fermi-Ulam theorems 
\par
\noindent
{\it Proof of theorem 2,1:}
We prove theorem 2,1 by applying the KAM theorem of Moser [M].
Namely let $(r,\theta)$ be polar coordinates of the annulus
$0 < a \le r \le b$ and difine a mapping of the annulus by
$$\eqalign{
\theta_1 &:= \theta + \alpha(r) + F(r,\theta)\cr
r_1 &:= r + G(r,\theta).\cr}
\eqno(4.1)$$
If $h(r,\theta)$ is a function with continuous derivatives up to order $s$
we define the $s$th derivative norm by
$$|h|_s := \sup \left | \left ( {{\del} \over {\del r}} \right )^{\sigma_1}
\left ( {{\del} \over {\del \theta}} \right )^{\sigma_2} h(r,\theta)
\right |, \quad \sigma_1 + \sigma_2 \le s\eqno(4.2)$$
where $(r,\theta)$ range over the domain in which $h$ is defined. 

\proclaim Theorem [M].  Fix $\e > 0$ and an integer $s \ge 1.$
Assume for the mapping (4.1) that every closed curve near a circle and its
image curve intersect.  Assume further $b-a \ge 1$ and
$$b_0^{-1} \le {{d\alpha(r)} \over {dr}} \le b_0\eqno(4.3)$$
with some constant $b_0 > 1.$  
Finally assume $F,G$ have continuous derivatives up to order $l$ and
satisfy the inequalities
$$\eqalign{
|F|_0 + |G|_0 &< \delta_0\cr
|\alpha|_l + |F|_l + |G|_l &< b_0\cr}
\eqno(4.4)$$
where $\delta_0 = \delta_0(\e,s,c_0)$ is a sufficiently small positive real 
number and $l = l(s)$ is a large enough integer. 
Then the mapping (4.1) has a closed invariant curve
$$\eqalign{
\theta &= \theta^{\prime} + p(\theta^{\prime})\cr
r &= r_0 + q(\theta^{\prime})\cr}
\eqno(4.5)$$
where the functions $p,q$ are functions of period $2\pi$ with $s$ continuous
derivatives satisfying
$$|p|_s + |q|_s < \e.\eqno(4.6)$$
\par

The application of Moser's theorem 
follows the same path as the proofs in [Pu4][Pu5] and 
[Pu9]. Let $(t,v)$ be polar coordinates for the plane
$\R^2$, $t \in
[0,T)$ the angular coordinate and $v \ge 0$ the radial coordinate. Let
$D_r \subset \R^2$ be the disk of radius $r$ and center $v=0.$  
For $r > 0$ sufficiently large let
$A: \R^2 \backslash D_r \r \R^2$ having the form $A(t,v) = (\pt,\pv)$ be 
given implicitly by the following formulas:
$$\eqalign{
f_2(\tt) - f_1(t) & = v \cdot (\tt - t) - {{g \cdot (\tt - t)^2} \over 2}\cr
\tv & = -v + g \cdot (\tt - t) + 2 \df_2(\tt)\cr
f_2(\tt) - f_1(\pt) & = - \tv \cdot (\pt - \tt) + {{g \cdot (\pt - \tt)^2} \over 2}\cr
\pv & = v + g \cdot (\pt - 2 \tt + t) + 2 \df_1(\pt) - 2 \df_2(\tt).\cr
}\eqno(4.7)$$
The formulas always have an implicit solution if $r$ is sufficiently large.
Here we take the minimal $\tt > t$ which is a solution of the first equation and
then the minimal $\pt > \tt$ which is a solution of the third equation.
Then the dynamics of the ball for large enough $v$ is implicitly given by the 
mapping $A.$ 
Namely, if at the moment of collision with the bottom plate the 
ball has large enough velocity then it will collide exactly once with the top 
plate before hitting the bottom plate again.
If these three events occur at times which we call $t < \tt < \pt$ and the
velocities at these times (at the moment after the collision) are called 
$v > 0, \ \tv < 0$ and $\pv > 0$ respectively then formula (4.7)
holds. 

We first verify that the intersection hypothesis of Moser's theorem holds.
\proclaim Lemma 4,1.
Let $\g$ be a simple closed smooth curve on $\R^2$ (parametrized by
$t$) enclosing the disk
$D_r$ for large enough radius $r.$ Then
$$\oint_{\g}\left ({{v^2}\over{2}} + f_1(t) g - v \df_1(t) \right ) dt
= \oint_{A\g}\left ({{v^2}\over{2}} + f_1(t) g - v \df_1(t) \right ) dt.
\eqno(4.8)$$
 
\noindent
{\it Proof of Lemma 4,1:}
We represent $A$ as the composition of four mappings $A = A_4 \circ A_3 \circ
A_2 \circ A_1$ where

$$\eqalign{
A_1(t,v) := & (t_1,v_1) = (\tt, v- g(\tt -t)),\cr
A_2(t_1,v_1) := & (t_2,v_2) = (\tt,\tv),\cr
A_3(t_2,v_2) := & (t_3,v_3) = (\pt,\tv - g(\pt-\tt)),\cr
A_4(t_3,v_3) := & (t_4,v_4) = (\pt,2 \df_1(\pt) - \tv + g (\pt - \tt) ).\cr
}\eqno(4.9)$$

The physical meaning of these mappings is as follows.
$A_1$ accounts for the movement from the lower plate to the upper plate,
$A_2$ accounts for the reflection due to the collision with 
the upper plate, $A_3$ accounts for the 
drop from the upper plate to the lower plate and $A_4$ for the reflection 
due to the collision with the lower plate.

The integral invariant of Poincar\'e-Cartan [C] gives

$$\eqalign{
&\oint_{\g} v dz - H(z,v) dt = - \oint_{A_1(\g)} v_1 dz_1 -
H(z_1,v_1) dt_1\cr
&-\oint_{A_2 \circ A_1 (\g)} v_2 dz_2 - H(z_2,v_2) dt_2 = 
- \oint_{A_3 \circ A_2 \circ A_1 (\g)} v_3 dz_3 - H(z_3,v_3) dt_3\cr
}\eqno(4.10)$$
where $z = f_1(t), H(z,v) = v^2/2 + gz, z_1 = z_2 = f_2(t_1) = f_2(t_2),
z_3 = f_1(t_3).$

Set $dz = \df_1(t), \ dz_1 = dz_2 = \df_2(t_1) dt_1 = \df_2(t_2) dt_2, \ 
dz_3 = \df_1(t_3) dt_3.$ This transforms equations (4.10) to

$$\eqalign{
\oint_{\g} \left ( v\df_1(t) - (v^2/2 + f_1(t))g \right )dt & = - \oint_{A(\g)}
\left ( v_1\df_2(t_1) - (v_1^2/2  + f_2(t_1)g) \right )  dt_1\cr
\oint_{A_2 \circ A_1(\g)} \left (  v_2 \df_2(t_2) - (v_2^2/2 + f_2(t_2)g) 
\right ) dt_2 & = \cr
= - \oint_{A_3 \circ A_2 \circ A_1(\g)}& \left ( v_3 \df_1(t_3) - 
(v_3^2/2 + f_1(t_3) g ) \right ) dt_3.\cr}
\eqno(4.11)$$

Equation (4.11) and $v_2 = 2 f_2(t_1) - v_1, \  t_1 = t_2, \  v_4 = 2 
\df_1(t_3) -v_3, \  t_4 = t_3$ yields

$$\oint_{\g} \left ( {{v^2} \over 2} f_1(t)g - v \df_1(t) \right ) dt = 
\oint_{A( \g)} \left ( {{(\pv)^2} \over 2} + f_1(\pt)g - v \df_1(\pt) \right ) 
d\pt.\eqno(4.12)$$

For large $r$ the the curves $A(\g)$  and $\g$ have the same orientation.
Thus the equality (4.12) is equivalent to the statement of Lemma 4,1.
\hfill\qed

\proclaim Lemma 4,2.
Let $\g$ be a parametrized simple closed smooth curve on $\R^2$ enclosing $D_r$ 
for large enough radius $r.$ Assume additionally 
that the mapping $A$ takes $\g$ to the curve $A(\g)$ of the same orientation.
Then $\g$ intersects $A(\g).$\par

Lemma 4,2 follows directly from Lemma 4,1. This completes the verification of 
the intersection hypothesis in Moser's theorem.  Next we must verify the 
small parameter condition (4.4) in Moser's theorem.  To do this we will
make several changes in coordinate as in [Pu4] and [Pu5]. Namely
let $l$ be a positive number. We introduce new polar coordinates
$(t,y) := U(t,v) := (t,2l/v)$ ($t$ is the angle variable, $y$ the radius 
variable).  We consider the mapping $A$ in these coordinates.  Equation
(4.7) transforms to $\pA = U \circ A \circ U^{-1}: (t,y) \r (\pt,\py)$

$$\py = {{2l} \over v}  
= y + \phi(t,y) := y - {{ y^2 (2 \df_1(\pt)-2 \df_2(\tt) + (\pt - 2 \tt + t
) g)}\over 
{2 l(t) + y(2 \df_1(\pt)-2 \df_2(\tt) + (\pt - 2 \tt + t) g)}}\eqno(4.13)$$
and there exists a $r > 0$ such that in the region $|y| \le r, 0 \le t \le
T$

$$\left |\phi(t,y)\right | < c_1y^2.\eqno(4.14)$$

Here $c_1$ is a positive constant not depending on $t$ or $y.$
Furthermore from (4.7) follows

$$\eqalign{
\tt -t = & {{2(f_2(\tt) - f_1(t))}\over{v + 
\sqrt{v^2 + 2 g (f_2(\tt) - f_1 (t))}}}\cr
\pt - \tt = & {{2(f_2(\tt) - f_1(\pt))}\over{- \tv + 
\sqrt{\tv^2 + 2 g (f_2(\tt) - f_1 (\pt))}}}; \quad ( \tv < 0).\cr
}\eqno(4.15)$$
Combining these two equations gives
$$\pt = t + {{f_2(\tt) - f_1(t)}\over {v}} + {{f_1(\pt) - f_2(\tt)}\over {\tv}}
+ K\eqno(4.16)$$

where 

$$K = {{- \alpha_1(f_2(\tt) - f_1(t))}\over {v(2v + \alpha_1)}} + 
{{\alpha_2(f_2(\tt) - f_1(\pt))}\over{\tv(2 \tv - \alpha_2)}}\eqno(4.17)$$ 

$$ \alpha_1 = {{K_1}\over{\sqrt{v^2 + K_1} + v}}, \ 
\alpha_2 = {{K_2}\over{\sqrt{\tv^2 + K_2} - \tv}}\eqno(4.18)$$

$$K_1 := -2g(f_2(\tt) - f_1(t)), \ K_2 := 2g(f_2(\tt) - f_1(\pt)).
\eqno(4.19)$$

Applying the change of variables $y= 2l/v$ to (4.7) and (4.16) yields

$$\tv = {{-2l + y ( g (\tt - t) + 2 \df_2(\tt))}\over{y}},\eqno(4.20)$$

$$\eqalign{
\pt = t + y + \psi(t,y) := t + y + \Bigl \{
{{y^2 l (g(\tt - t) + 2 \df_2(\tt))}\over{2l - y (g (\tt - t)  + \df_2(\tt))}} 
& + {{y(f_2(\tt) - l - f_1(t))}\over{2l}}+ \cr
 + & {{y(f_2(\tt) - l - f_1(\pt))}\over{2l-y ( g (\tt - t) + 2\df_2(\tt))}} 
+ K
\Bigr \}.\cr}\eqno(4.21)$$

Next choose $l$ large enough so that

$$l > \sup_{t_1,t_2} |f_2(t_2) - f_1(t_1)|.\eqno(4.22)$$

Then for sufficiently small $r > 0$ 
$$\left | \psi(t,y)\right | < c_0 |y|, \quad 
\left | {{ \del \psi} \over {\del y}}
\right | < c_0, \eqno(4.23)$$
for $|y| < r, \ 0 \le t < T$ where $0 < c_0 < 1$ is a constant not depending 
on  $t,y.$

If the constants
$c_0$ and $c_1$ in equations (4.23) and (4.14) where sufficiently
small then Moser's theorem would imply that in any neighborhood
of the point $y=0$ there exists an $A^{\prime}$-invariant curve
which surrounds the point $y=0.$ Theorem 2,1 follows from this fact.

In our case the constants $c_0$ and $c_1$ from (4.23) and (4.14)
are not necessarily small.  The only inequality which $c_0$
satisfies is $0 < c_0 < 1.$  Thus we apply the proof of theorem
2 in [Pu5] which is a generalization of Moser's theorem in the
case that the mapping has the form given by equations (4.13) and 
(4.21) which satisfy Moser's intersection condition and satisfy
the inequalities (4.14) and (4.23) with an arbitrary constant $c_1$
and $0 < c_0 < 1.$  The main part of the proof of theorem 2 in
[Pu5] is the construction of change of variables which transforms
the mapping $A^{\prime}$  into a mapping which satisfies 
Moser's theorem, that is the newly constructed mapping satisfies
Moser's intersection condition of has constants analogous
to $c_0$ and $c_1$ which are small. 
This fact finishes the 
proof of theorem 2,1.\hfill\qed

\noindent
{\it Proof of Theorem 2,2:}
Equation (4.7) implies
$$\eqalign{
\tt &= t + {{f_2(\tt) - f_1(t)}\over{v(1 - {{g}\over{2v}} 
(\tt - t))}}\cr
\pt &= \tt + {{f_2(\tt) - f_1(\pt)}\over{-\tv(1 - {{g}\over{2\tv}} 
(\pt - \tt))}}\cr
}\eqno(4.24)$$

We assume without loss of generality that $T = 1, \ t^* = 0$ and
$f_1(0) = f_2(0) = 0.$ The assumption of the theorem imply that
$f^{(k-1)}_1(0) =  f^{(k-1)}_2(0)$ and 
$f^{(k)}_1(0) <  f^{(k)}_2(0)$ for $k=1$ or 2. 

Choose $\epsilon > 0$ so small that
$$\left | \sum_{j=0}^k(f^{(j)}_i(0)/j!)t^j - f_i(t) \right | 
< 2 \max ( 1, |f_1^{(k+1)}(0)|, 
|f_2^{(k+1)}(0)|) t^{k+1}\eqno(4.25)$$
for $ t \in (-\epsilon,0], 
\ i=1,2.$ Fix a positive constant $C_1.$  
We choose the constant $K$ so large that if $(v_0,t_0)$ satisfies $|v_0| > K$ 
and $t_0 \in (-\epsilon, 0)$ then
the ball will in its three next hits alternate hitting the top and bottom 
plates.
We denote the next hit of the bottom plate by $(v,t)$ (the velocity is taken
to be the velocity just after the hit) and the next hit after time 
$t$ of the top 
plate by $(\tv,\tt)$ (again the velocity is that just after the hit). We 
further assume $K$ is so large that $v,-\tv > 2C_1.$ 

If $ -\epsilon < t < 0$ then
$$\eqalign{
{-\tv \left ( 1 - {{g}\over{2\tv}} (\pt - \tt) \right ) } > C_1\cr
{v \left ( 1 - {{g}\over{2v}} (\tt - t) \right ) } > C_1.\cr
}\eqno(4.26)$$

Putting together (4.24) and (4.26) yields
$$\pt \le t + {{\left | f_2(\tt) - f_1(t)\right |}\over{C_1}} + 
{{\left | f_2(\tt) - f_1(\pt) \right |}\over{C_1}}.\eqno(4.27)$$
Now equations (4.24) and (4.26) imply $\max(|\tt^k - t^k|,|\tt^k-{\pt}^k|) = 
O(t^{k+1})$ and thus
$$\max(f_2(\tt) - f_1(t),f_2(\tt) - f_1(\pt)) = 
{{f^{(k)}_2(0) - f^{(k)}_1(0}) \over k!} t^k + O(t^{k+1}).
\eqno(4.28)$$
Putting (4.27) and (4.28) together yields

$$\pt < t + {{f^{(k)}_2(0) - f^{(k)}_1(0)} \over {k!C_1}} t^k + O(t^{k+1}).
\eqno(4.29)$$

We inductively assume that $v_i > K$ for $i = 1,2,\dots n-1.$
Let $(t_n,v_n) = A^n(t,v)$, that is the time and velocity of the $n$th 
hit on the bottom plate.  Furthermore let $\tt_n$ play the role of $\tt$ 
when $t = t_n$ in equation (4.7). Equations (4.24) and (4.26) imply that
$\max(|t_n-\tt_n|,|t_{n+1}-\tt_n|) = O(t_n^2)$ and (4.7) implies
$$\eqalign{
v_{n+1} &= v_n + 2\df_1(t_{n+1}) - 2 \df_2(\tt_n) + O(t_n^2)\cr
& =  v_n + 2(f^{(k)}_1(0)t^{k-1}_{n+1} - f^{(k)}_2(0)\tt^{k-1}_n) + O(t_n^k)\cr
& = v_n + 2 (f^{(k)}_1(0) - f^{(k)}_2(0))t^{k-1}_n + O(t_n^k)\cr
& = v_0 + 2 (f^{(k)}_1(0) - f^{(k)}_2(0)) \sum_{i=1}^n (t^{k-1}_i + 
O(t_i^k))\cr
& \ge v_0 + 2 (f^{(k)}_1(0) - f^{(k)}_2(0)) \sum_{i=1}^n (t^{k-1}_i + 
(1/2)|t^{k-1}_i|).\cr }
\eqno(4.30)$$

Equation (4.30) implies that $v_{n} \ge v_{n-1} > K,$ verifying the inductive
assumption.  
If $k=1$ then it is clear that $v_n \r \infty.$ If $k=2$ then 
let $s_0 := t_0$ and $s_{n+1} := s_n + {{\ddf_2(0) - 
\ddf_1(0)} \over {3C_1}} s_n^2.$
Then the following lemma implies that $\sum s_n = - \infty.$ 
If $-t_0 > 0$ is sufficiently small then equation (4.29) implies that
$t_{n+1} < t_n + {{\ddf_2(0) - \ddf_1(0)} \over {3C_1}} t_n^2$
and then simple algebra yields that
$t_n \le s_n$ for all $n.$ Thus it follows that $\sum t_n = -\infty$ and 
$v_n \rightarrow \infty.$
\hfill\qed

\proclaim Lemma 4,3.
Supposed $\tau(s) \not \equiv 0$ is a positive $C^3$ function in the domain 
$(-\delta,0]$ where $\delta$ is a positive constant.  Furthermore we assume 
$\tau(0) = {\dot \tau}(0) = 0.$ Consider the sequence 
$\{s_n: n \ge 0\}$ 
given by $s_{n+1} = s_n + \tau(s_n).$  Then there exists an $\epsilon_0 > 0$
so that if $-s_0 \in (0,\epsilon_0)$ then $\sum_{n=0}^{\infty} (-s_n)
= \infty.$\par

\noindent
{\it Proof of Lemma 4,3:}
Suppose that 
$$\sum_{n=0}^{\infty} (-s_n) < \infty.\eqno(4.31)$$ 
\noindent
Then 
$$\sum_{k=0}^{\infty} {{-\tau(s_k)}\over{s_k}} < \infty \eqno(4.32)$$
since $\tau(0) = {\dot \tau}(0) = 0.$ For any $n$ we have the identity
$$s_{n+1} = s_0 \prod_{k=0}^n \left ( 1 + {{\tau(s_k)}\over {s_k}} 
\right ).\eqno(4.33)$$

Thus by taking logarithms of (4.32) we see that the limit
$$\lim_{n \rightarrow \infty} \prod_{k=0}^n 
\left ( 1 + {{\tau(s_k)}\over{s_k}} \right ) = 
s^*\eqno(4.34)$$
exists and $s^* \ne 0.$ Therefore, from (4.33)
$\lim_{n \r \infty} s_n = s_0 s^* \ne 0$
which contradicts equation (4.31).\hfill\qed

\section 5. Proofs of bouncing ball theorems\par
\noindent
{\it Proof of proposition 3,1:}
Suppose the masses of the balls are $m_1$ and $m_2.$
Let $$\alpha := {{m_1 - m_2}\over {m_1 + m_2}}.\eqno(5.1)$$
Let the velocities of the balls just before a collision be given by
$\vo_-,\v2_-$ and just after the collision by $\vo_+,\v2_+.$
Then the following equations hold
$$\eqalign{
\vo_+ = & \alpha \vo_- + (1 - \alpha) \v2_-\cr
\v2_+ = & (1 + \alpha) \vo_- + - \alpha \v2_-.\cr}\eqno(5.2)$$
Then proposion follows from equation (5.2).
\hfill\qed

\noindent
{\it Proof of theorem 3,2:}
We first assume there is only one ball above the plate and will build
two trajectories $\ga2_1,\ga2_2$ for it, where $\ga2_1$ is periodic and $\ga2_2$
has unbounded velocity.
Suppose that after the collision with the plate at time $t$ the ball 
has velocity $v > 0,$ after the next collision with the plate at 
time $\pt$ the ball has velocity $\pv$ and in the time interval
$[t,\pt]$ the ball does not collide with the plate. 
The dynamics of the ball is implicitly
given by the mapping $A:(t,v) \r (\pt,\pv)$:

$$\eqalign{
f(t) + v \cdot (\pt - t) - {g \over 2} \cdot
(\pt - t)^2 = & f(\pt)\cr
g \cdot (\pt - t) - v + 2 \df(\pt)=&\pv.\cr
}\eqno(5.3)$$
Here we take the minimal $\pt > t$ which is a solution of the first
equation. If $f \in C^1$ then there exists $\bt2_0$ so that $\df(\bt2_0) = 0.$
If $f \in \C$ then there exists $\bbt2_0$ so that $\df(\bbt2_0) = Tg/2.$ 
Set $v_0 := Tgm_2/2.$ Here we have assumed that $K=1$ in the definition of the
class $\C.$

Using the points $\bt2_0, \bbt2_0$ and elementary algebra we produce 
two trajectories, 
$$\eqalign{
\ga2_1 & := \{ (\bt2_0,v_0), (\bt2_1,\bv2_1), \cdots, 
(\bt2_n,\bv2_n), \cdots \},\cr
\ga2_2 & := \{ (\bbt2_0,v_0), (\bbt2_1,\bbv2_1), \cdots, 
(\bbt2_n,\bbv2_n), \cdots \}\cr
}\eqno(5.4)$$
such that $v_0 = Tgm_2 /2 > 0$ where $m_2$ is a sufficiently large
positive integer,
$\bt2_0, \bbt2_0$ satisfy the equations 
$\df(\bt2_0) = 0, \ \df(\bbt2_0) = Tg/2$ and for all $n \ge 0$
the following hold:
$$\eqalign{
\bv2_{n+1} = v_0, \ & \bbv2_{n+1} = \bbv2_n + Tg,\cr
t_0 = &\bt2_{n+1} = \bbt2_{n+1} \hbox{ mod}(T),\cr
\bt2_{n+1} - \bt2_n = Tm_2, \ &
\bbt2_{n+1} - \bbt2_n = T(m_2 + 2n).\cr
}\eqno(5.5)$$

Now we return to the case of two balls.  The case
when $\lim_{i \r \infty} |\vo_i| = \lim_{i \r \infty} |\v2_i| = \infty,$
occurs when the initial condition of ball $P_1$ is $(\bbt2_0,\bbv2_0)$
and the initial condition of ball $P_2$ is 
$(\bbt2_0+nT,\bbv2_0)$ for any positive integer $n.$ 
The case $\limsup_{i \r \infty} |\vo_i| =
\limsup_{i \r \infty} |\v2_i| = \infty,$
but $\max (\liminf_{i \r \infty} |\vo_i|,
\liminf_{i \r \infty} |\v2_i|) < \infty,$
occurs when the initial condition of the balls is
$(\bbt2_0,\bbv2_0)$ and $(\bt2_0,\bv2_0).$  In both cases the behavior
described occurs because when two ball of equal mass collide they exchange
velocities.
\hfill\qed

\noindent
{\it Proof of theorem 3,3:}
Let $\vo(t),\v2(t)$ be the velocities of balls $P_1,P_2$ at time
$t.$ Suppose that there is a constant $C_1 > 0$ such that
$\sup_{t \ge 0} |\vo(t)| < C_1$ and $\sup_{t \ge 0} |\v2(t)| = \infty.$
Let $0 \le t_1 < t_2 < \dots$ be the collision times between the particles
and $\Delta t_i = t_{i+1} - t_i.$

If $\sup_{i \ge 0} \Delta t_i = \infty$ then there is an $i \ge 0$ for which
$\Delta t_i > 4C_1 / g.$  Then for $t_* = t_i + (1/2) \Delta t_i$ we have
$$\vo(t_*) = \vo(t_i) - g{{\Delta t_i}\over 2}
\le C_1 - 2C_1 < -C_1\eqno(5.6)$$
a contradiction.

It remains to study the case when $\exists C_2 > 0$ such that
$\sup_{i \ge 0} \Delta t_i < C_2.$ 
Since the second ball has unbounded velocity there are arbitrarily large
time intervals in which it has arbitrarily large speed.  In particular there
exists $\hat t$ such that for all $t \in [\hat t,\hat t + C_2]$
we have 
$$|\v2(t)| > C_1.\eqno(5.7)$$  
But since $\Delta t_i < C_2$ there is always a $t_i$ in the interval
$[\hat t,\hat t + C_2].$ 
Now $\vo(t_i) = -\vo(t_i^-) + 2 \v2(t_1).$ Thus using (5.6) and
(5.7) yields
$$|\vo(t_i)| > -C_1 + 2C_1 = C_1\eqno(5.8)$$
a contradiction.

We have already shown in theorem 3,2 that for a single ball
there is an unbounded
trajectory $\ga2_2$ if $f \in \C.$  However the
zero mass particle does not effect the dynamics of the ball $P_2$
so this trajectory occurs in the case (1) as well. The above
discussion shows that any initial condition for $P_1$ above the
trajectory $\ga2_2$ 
leads to both balls having unbounded velocity.
\hfill\qed

\noindent
{\it Proof of theorem 3,4}
We apply Theorem 2,2. The role of the function $f_2(t)$ in theorem 2,2 will
be played by the orbit of $P_2$ and the role of 
$f_1(t)$ will be played by $f(t).$  
\hfill\qed

\section 6. Acknowledgements
\par
LDP thanks BiBoS Research Center for the invitation to visit Bielefeld.
SET is grateful for the support of the Deutsche Forschungsgemeinschaft.

\section 7. References
\par
\frenchspacing
\baselineskip=16pt

\item{[A]} V.I.~Arnold {\it On the behavior of an adiabatic invariant
under slow periodic variation of the Hamiltonian} Sov.~Math.~Dokl.~{\bf 3}
(1962) pp 136-140. 

\item{[C]} \'E.~Cartan {\it Lecons sur les invariants int\'egraux}
Hermann, Paris (1922).

\item{[Ch]} N.~Chernov {\it The ergodicity of a Hamiltonian system of two
particles in an external field} Physica D {\bf 53} (1991) pp 233-239.

\item{[F]} E.~Fermi {\it On the origin of the cosmic radiation} 
Phys.~Rev.~{\bf 75} (1949) pp 1169-1174.

\item{[M]} J.~Moser {\it On invariant curves of area preserving mappings of
an annulus} Nachr.~Akad. Wiss.~G\"ottingen II Math.-Phys.~{\bf 1} (1962)
pp 1-20.

\item{[P]} H.~Poincar\'e {\it R\'eflexions sur la th\'eorie cin\'etique des 
gaz} J.~Phys.~th\'eoret.~appl.~{\bf 4} ser.~5 (1906) pp 369-403.

\item{[Pu1]} L.D.~Pustyl'nikov {\it On a certain dynamical system}
Russ.~Math.~Surv.~{\bf 23} (1968).

\item{[Pu2]} L.D.~Pustyl'nikov {\it Existence of a set of positive measure 
of oscillating motions in a certain problem of dynamics} 
Sov.~Math.~Dokl.~{\bf 13} (1972) pp 38-40.

\item{[Pu3]} L.D.~Pustyl'nikov {\it Stable and oscillating motions in 
nonautonomous dynamical systems II} Mosc.~Math.~Soc.~{\bf 34} (1977)
pp 1-101.

\item{[Pu4]} L.D.~Pustyl'nikov {\it On the Ulam problem} 
Theor.~Math.~Phys.~{\bf 57} (1983) pp 1035-1038.

\item{[Pu5]} L.D.~Pustyl'nikov {\it On the Fermi-Ulam model} 
Sov.~Math.~Dokl.~{\bf 292} (1987) pp 88-92.

\item{[Pu6]} L.D.~Pustyl'nikov {\it A new mechanism of particle 
acceleration and a relativistic analog of the Fermi-Ulam model}
Theor.~Math.~Phys.~{\bf 77} (1988) pp 1110-1115.

\item{[Pu7]} L.D.~Pustyl'nikov {\it A new mechanism of particle 
acceleration and rotation numbers} Theor.~Math.~Phys.~{\bf 82} (1990)
pp 180-187.

\item{[Pu8]} L.D.~Pustyl'nikov {\it A mechanism if irreversibility and 
unbounded growth of the energy in a model of statistical mechanics}
Theor.~Math.~Phys.~{\bf 86} (1991) pp 82-89.

\item{[Pu9]} L.D.~Pustyl'nikov {\it On the existence of invariant curves 
for maps close to degenerate and a solution of a problem of Fermi-Ulam}
Math.~Sbor.~(1994) pp 1-12 (in Russian).

\item{[U]} S.~Ulam {\it On some statistical properties of dynamical systems}
Proc. Fourth Berkeley Symp.~Math.~Stat.~Prob.~vol.~3, Univ.~Cal.~Press,
Berkeley, CA.~(1961).

\item{[W1]} M.~Wojtkowski {\it A system of one dimensional balls with gravity} 
Comm.~Math.~Phys. {\bf 126} (1990) pp 507-533.

\item{[W2]} M.~Wojtkowski {\it A system of one dimensional balls in an external
field II} Comm.~Math. Phys.~{\bf 127} (1990) pp 425-432.

\vfill\eject
\centerline{Captions for figures.}
\item{Figure 1} The Fermi-Ulam model in an external gravitational field.
\item{Figure 2} Touching plates.
\item{Figure 3} Zero mass ball on top.
\item{Figure 4} Zero mass ball in the middle.
\end